\documentclass[12pt,a4paper]{amsart}

\usepackage{amsmath, amsfonts, xifthen, latexsym, amssymb, amsthm, amscd}
\usepackage[utf8]{inputenc}
\usepackage{graphicx}
\usepackage{url}

\usepackage{hyperref}
\hypersetup{colorlinks=true,citecolor=blue,filecolor=blue,linkcolor=blue,urlcolor=blue}
\usepackage[margin=1.4in]{geometry}

\newtheorem{theorem}{Theorem}
\newtheorem{lemma}{Lemma}[section]

\newtheorem{proposition}{Proposition}[section]

\newtheorem{corollary}{Corollary}[section]

\theoremstyle{definition}
\newtheorem{definition}{Definition}[section]

\newcommand{\eps}{{\varepsilon}}

\renewcommand{\phi}{{\varphi}}
\newcommand{\grd}{{Gr_d}}

\newcommand{\Supp}{\operatorname{Supp}}

\newcommand\Sep{\operatorname{Sep}}

\newcommand{\cB}{\mathcal{B}}\newcommand{\cP}{\mathcal{P}}\newcommand{\cD}{\mathcal{D}}

\newcommand\R{\mathbb R}

\newcommand{\cA}{\mathcal{A}}
\newcommand{\cQ}{\mathcal{Q}}
\newcommand{\cS}{\mathcal{S}}
\newcommand{\cV}{\mathcal{V}}
\newcommand{\cU}{\mathcal{U}}

\newcommand{\cC}{\mathcal{C}} \newcommand{\cL}{\mathcal{L}}
\newcommand{\cH}{\mathcal{H}}

\newcommand{\vi}{\vskip 0.1in \noindent}
\usepackage{etoolbox}
\newtoggle{no_cases}\toggletrue{no_cases}
\newcommand{\case}[2][]{\iftoggle{no_cases}{\left\{\begin{array}{ll}#2 & #1}{\\#2 & #1}\togglefalse{no_cases}}
\newcommand{\esac}{\end{array}\right.\toggletrue{no_cases}}

\newcommand{\Prob}{\operatorname{Prob}}

\begin{document}
\title[Approximate Proof Labeling Scheme for Planarity]{Planarity can be Verified by an Approximate Proof Labeling Scheme in Constant-Time} 
\author{Gábor Elek}
\address{\vi Department of Mathematics And Statistics, Fylde College, Lancaster University, Lancaster, LA1 4YF, United Kingdom   
\vskip 0.05in \noindent and
\vskip 0.05in \noindent Alfred Renyi Institute, Budapest, Hungary}

\email{g.elek@lancaster.ac.uk}  

\thanks{The author was partially supported
by the ERC Consolidator Grant "Asymptotic invariants of discrete groups,
sparse graphs and locally symmetric spaces" No. 648017. }

\begin{abstract}
Approximate proof labeling schemes  were introduced by \\Censor-Hillel, Paz and Perry \cite{CPP}. Roughly speaking,
a graph property~$\cP$ can be verified by an approximate proof labeling scheme in constant-time if the vertices of a graph having the property can be convinced, in a short period of
time not depending on the size of the graph, that they are having the property $\cP$ or at least
they are not far from being having the property $\cP$. The main result of this paper is that bounded-degree planar graphs (and also outer-planar graphs, bounded genus graphs, knotlessly embeddable graphs etc.) can be verified by an approximate proof labeling scheme in constant-time.

\end{abstract}\maketitle
\newpage

\setcounter{tocdepth}{2}
\tableofcontents
\vi
\vi
\section{Introduction}
\vi
Our paper is about constant-time \textbf{distributed graph algorithms} introduced in the
seminal work of Naor and Stockmeyer \cite{Naor}. Such an algorithm runs on simple undirected
graphs $G$ of bounded degree. Each vertex of $G$
collects information from vertices located within radius $r$ from it and creates
some output based on the collected information. In a \textbf{distributed decision algorithm}
each vertex outputs a decision: accept or reject. Then, the graph $G$ is accepted by
the collective decision (verification) of the vertices if all the vertices accept and
the graph is rejected if at least one of the vertices rejects.
Unfortunately, there are not too much interesting graph families that
can be verified in this way. However, one can consider a non-deterministic version
of graph verification: \textbf{proof labeling schemes} (introduced by Korman, Kutten and Peleg \cite{KKP})
or under the name of \textbf{locally checkable proofs} (due to
G\"o\"os and Suomela \cite{GS}) with somewhat different
conditions (see \cite{Feuil} for an
extended survey).
Here, a prover helps the vertices to make their decision.
\vi
The prover labels the vertices with an element of a finite set $Q$ and the
vertices can view their $r$-neighbourhoods as $Q$-labeled balls for some positive integer $r$. A graph family (note that in
our paper we only consider graph families of bounded degrees) can
be verified by a proof labeling scheme if there exists
a labeling-verification protocol such that
\vi
$\bullet$ If a graph $G$ is in the family then there exists
a labeling that makes all vertices accept.
\vi
$\bullet$ If a graph $G$ is not in the family then for
all labelings there exists at least one vertices that rejects.
\vi
Proof labeling schemes help to verify more interesting graph classes. Nevertheless,
one the most interesting classes, namely, the class of planar graphs still cannot
be verified by a proof labeling scheme in constant-time \cite{FFP}. Our graphs are of bounded-degree, the number of labelings
are finite and the nodes explores a constant distance neighbourhood of themselves, that is why we use the \emph{constant-time} terminology.  \vi
However, a relaxation of the above verification procedure, approximate proof labeling schemes,
was introduced by Censor-Hillel, Paz and Perry \cite{CPP}(see also \cite{EG}).
In the case of an approximate proof labeling scheme the vertices may accept a graph
even if it is not in the graph family, provided that it is not far from the family in edit distance. One should note
that in \cite{CPP} the proof labeling schemes are used to certificate the approximation of
certain parameters such as the maximum matching or a maximum independent set, nevertheless the idea is very similar. 
\vi
\textbf{Our result  }   We will show (Theorem \ref{main}) that the class of \textbf{planar graphs} and
in general all \textbf{monotone hyperfinite} (see Section \ref{hyper}) graph classes such as outer-planar
graphs, bounded genus graphs or  knotlessly embeddable graphs (or any other\textbf{ minor-closed families}) can be verified
with an approximate proof labeling scheme in constant-time. 
\vi
\textbf{Related work}  It is important to note that the class of planar graphs can be verified 
with a proof labeling scheme if we allow the vertices to be labeled by $O(\log(n))$-bits strings,
where $n$ is the size of the graph \cite{FFP}, even without the
bounded degree assumption. This result can be extended to the class of bounded genus graphs
as well (\cite{FFP2} and \cite{EL}).
\vi
Recently, Romero, Wrochna and \v{Z}ivn\'y \cite{RWZ} constructed polynomial-time approximation schemes for certain maximum constraint
satisfaction problems in the case of monotone hyperfinite graph classes. The main novelty
of their approach was the application of
strong hyperfiniteness, a strengthening of the
hyperfiniteness property. Our proof is
also based on strong hyperfiniteness, but in the
form of \textbf{Property A.} This geometric
property can be used for
proof verification in a natural way and had
already important applications in group theory,
in algebraic topology and in the theory of operator
algebras. 

\section{Approximate proof labeling schemes} \label{masodiksection}
\vi In order to avoid any confusion, let us fix some terminologies.
For an integer $d>1$, let $\grd$ be the set of
all finite, simple graphs of maximum degree at most $d$. If $x,y$ are adjacent
vertices in a graph $G$, we use the notation $x\sim y$.
For all graphs $G\in\grd$ we will consider the shortest
path distance $d_G$ on the vertex set $V(G)$ of $G$.
In our paper we consider only properties $\cP$ such
that $\cP\subset \grd$ for some $d$.
\vi
By a \textbf{ball} of radius $r$, we mean a
finite connected graph $B$ with a distinguished
vertex $v$ (the center) such that
$$\max_{y\in V(G)} d_G(v,y)=r\,.$$
\vi
For a fixed graph $G$ and vertex $x\in V(G)$ the 
\textbf{neighbourhood of radius $s$} centered at $x$
is the subgraph $B_s(x,G)$ induced on the vertices $y$ such that
$$d_G(x,y)\leq s\,.$$
\vi
It is important to note that a neighbourhood $B_s(x,G)$ 
above is always a ball with center $x$, however,
the radius of this ball can be equal to $s$ only if
the diameter of the graph $G$ is at least $s$. If the
diameter of $G$ is less than $s$, then the radius of $B_s(x,G)$
as a ball is always less than $s$. We will denote
by $N^d_r$ the maximum size of a ball of radius $r$
with maximum degree at most $d$.
\vi
For a graph $G\in\grd$ and a finite set $Q$, a \textbf{$Q$-proof}
is a function $T:V(G)\to Q$.\vi
A \textbf{$Q$-verifier} $\cV$ of local horizon $r$ is a subset of $B^{Q}_{r,d}$,
where $B^{Q}_{r,d}$ is the set of all $Q$-vertex labeled balls of radius at most $r$ and maximum degree at most $d$.
\vi
A $Q$-verifier $\cV$ of local horizon $r$ \textbf{accepts} 
a $Q$-proof $T$ on the graph $G\in\grd$, if
for all vertices $x\in V(G)$, $B_r(x,G,T)\in \cV$,
where $B_r(x,G,T)$ is the  neighbourhood of radius $r$
centered at $x$ with vertex labelling induced by $T$.\vi
A $Q$-verifier $\cV$ of local horizon $r$ \textbf{rejects} 
a $Q$-proof $T$ on the graph $G\in\grd$, if
for at least one vertex $x\in V(G)$, $B_r(x,G,T)\notin \cV$.\vi
We refer to subsets of $\grd$ as ``properties" and
we say that a property $\cP\subset\grd$ can be verified by a
\textbf{proof labeling scheme} (PLS) in constant-time if there exists a finite
set $Q$, a positive integer $r$
and a verifier $\cV\subset B^{Q}_{r,d}$
such that 
\begin{itemize}
\item for any $G\in\cP$ there exists a $Q$-proof
$T:V(G)\to Q$ accepted by $\cV$, 
\item for any $H\notin \cP$ all the $Q$-proofs on $H$
are rejected by $\cV$. 
\end{itemize} 
\vi
Verifiability by a PLS in constant-time entails that
the vertices of a graph $G \in \cP$ can be convinced
in a short period of time, that they are indeed vertices
of a graph having the given property. Clearly, $3$-colorability is such a property. Indeed, the proof $T:V(G)\to \{a,b,c\}$
will be the $3$-coloring and the verifier will check the
properness of the coloring on balls of radius $1$.
On the other hand,  Feuilloley et al. showed (\cite{FFP}, Theorem 2.) that \textbf{planarity cannot be verified by a PLS in constant-time.} 
\vi
In light of the result above we need a relaxation of the proof labeling
scheme verification procedure. 
Such relaxation has been introduced by Censor-Hillel, Paz and Perry \cite{CPP} (see also \cite{EG}) under the name approximate proof labeling scheme. \vi
First we need some terminology.
Recall that if $\cP\subset \grd$ and $H\in \grd$, then
the edit distance between the monotone property $\cP$ and the
graph $H$ is defined by
$$e(H,\cP):=\inf_{G\in \cP\,, V(G)=V(H)} 
\frac{|E(G)\triangle E(H)|}{|V(H)|}\,.$$
\noindent
\begin{definition}[\textbf{Approximate Proof Labeling Scheme}] A property $\cP\in 
\grd$ can be verified by an approximate proof labeling scheme in constant-time if for any $\eps>0$  there exists a set $Q_\eps\geq 0$, some positive constant $r_\eps$ and a verifier $\cV_\eps\subset B^{Q_\eps}_{r_\eps,d}$
such that
\begin{itemize}
\item for any $G\in\cP$ there exists a $Q_\eps$-proof $T$ on $G$
accepted by $\cV_\eps$,
\item for any $H, e(H,\cP)>\eps$,  all $Q_\eps$-proofs on
$H$ are rejected by $\cV_\eps$.
\end{itemize}
\end{definition} 
The main result of this paper is the following theorem.
\begin{theorem} \label{main}
Planarity and, in general, all monotone hyperfinite properties that
are closed under taking disjoint unions can be verified by an approximate proof labeling scheme in constant-time for bounded-degree graphs. 

\end{theorem}
\vi Note that a property is monotone if it is
closed under taking subgraphs. Hyperfiniteness will
be discussed in Section \ref{hyper}.
 \vi In \cite{BSS} the authors showed that
every minor-closed property is monotone hyperfinite (Theorem 1.1). In their paper they list several minor-closed properties
e.g. planarity, outer-planarity, graphs with bounded genus or
bounded tree-width.  By definition, all of these
properties are closed under taking disjoint unions. Hence, by our main theorem all of these
properties can be verified  with an approximate
proof labeling scheme in constant-time. 
\vi
Now, we introduce the notion of relative verifiability by PLS.
\begin{definition}
Let $\cP\subset \mathcal{Q}\subset \grd$ be properties.
We say that $\mathcal{P}$ can be verified by a PLS with respect to $\cQ$ in constant-time
if there exists a verifier $\cV\subset
B^{Q}_{r,d}$ such that
\begin{itemize}
\item for every $G\in\cP$ there exists a proof $T:V(G)\to Q$
such that $\cV$ accepts $T$,
\item for every $H\notin\mathcal{Q}$ all proofs $T:V(H)\to Q$ are rejected by $\cV$.
\end{itemize}
\end{definition} 
\vi
Now, let $G\in\grd, r\geq 1,$ $Q$ be a finite set and $\cV\subset 
B^{Q}_{r,d}$. We say that $\cV$ \textbf{verifies} $G$
if there exists a proof $T:V(G)\to Q$ such
that $\cV$ accepts $T$.
\vi
Let $\cV\subset B^{Q}_{r,d}$ be a verifier. We denote by $\cL_\cV$ the
set of graphs in $\grd$ verified by $\cV$. So, a property $\cP$ can be verified by a PLS in constant-time if there exists $\cV$ such that
$\cP=\cL_\cV$. \begin{lemma} \label{szorzat} If $\cV_1\subset B^{Q_1}_{r_1,d}$
and $\cV_2\subset B^{Q_2}_{r_2,d}$, then there exists
a verifier $\cV_3\in B^{Q_1\times Q_2}_{r_3,d}$ such that
$\cL_{\cV_3}=\cL_{\cV_1}\cap \cL_{\cV_2}$ and $r_3=\max(r_1,r_2)$. \end{lemma}
\proof
Let $B$ be a $Q_1\times Q_2$-labeled ball of radius
at most $r_3$. Let $B_1$ be the $r_1$-neighbourhood
of the center of $B$
equipped with the $Q_1$-labeling inherited from the
$Q_1\times Q_2$-labeling of $B$
using the first coordinate. Similarly,
we can define $B_2$. Let $B\in \cV_3$ if
$B_1\in\cV_1$ and $B_2\in\cV_2$.
\vi 
Assume that $G\in \cL_{\cV_3}$ and the proof
$T=T_1\times T_2:V(G)\to Q_1\times Q_2$
is accepted by $\cV_3$.
Then by definition, $T_1:V(G)\to Q_1$ is accepted
by $\cV_1$ and $T_2:V(G)\to Q_2$ is accepted
by $\cV_2$. Hence, $G\in \cL_{\cV_1}\cap \cL_{\cV_2}$.
\vi
Conversely, let $G\in \cL_{\cV_1}\cap \cL_{\cV_2}$
and let $S_1:V(G)\to Q_1$ be accepted
by $\cV_1$ and $S_2:V(G)\to Q_2$ be accepted
by $\cV_2$. Then by definition,
$$S=S_1\times S_2:V(G)\to Q_1\times Q_2$$
is accepted by $\cV_3$, thus $G\in \cL_{\cV_3}\,.$ \qed
\vi
By definition, if $\cP \subset\cQ$, then there exists $\cV$ such that
$$\cP\subset \cL_\cV \subset \cQ$$
\noindent
if and only if $\cP$ can be verified by a PLS in constant-time relative to $\cQ$.
If $\cP$ is a property let $\cP_\eps$ be the set of
graphs which are at most $\eps$-far in edit-distance from having the property $\cP$. Then $\cP$ can be verified by an approximate proof labeling scheme in
constant-time if for any $\eps>0$ there exists 
a verifier $\cV_\eps$ such that $\cP\subset \cL_{\cV_\eps} \subset \cP_\eps$.
\vi
In the introduction we mentioned a second verification protocol introduced by G\"o\"os and Suomela \cite{GS} under the name of locally checkable proofs. Sometimes a unique identifier
is provided for each node of the graph and the verifier
can take the identifiers into consideration as well. Although there are
certain properties for which one can benefit from the
existence of the identifiers, it follows from (\cite[Theorem 1.]{FHK}) that if one cannot verify a monotone property
without identifiers in constant-time, then one cannot verify
that property in constant-time even if unique identifiers are provided. 
\section{Hyperfiniteness} \label{hyper}
\vi
First, we recall the notion of hyperfiniteness (see \cite{BSS}) that
plays an important role in our paper.\vi
For $\eps>0$ and $K\geq 1$, a graph $G\in\grd$ is called
$(\eps,K)$-hyperfinite if there exists $W\subset E(G)$
such that
\begin{itemize}
\item $|W|\leq \eps |E(G)|$,
\item if we remove $W$ 
from $G$, in the remaining graph all the components
have size at most $K$.
\end{itemize}
\vi
A property $\cP\subset\grd$ is $(\eps,K)$-hyperfinite 
if all $G\in\cP$ are $(\eps,K)$-hyperfinite. The set
of all $(\eps,K)$-hyperfinite graphs of maximum degree at most $d$ is denoted by
$\cH^d_{\eps,K}$.
\vi
We call a property $\cP\subset\grd$ \textbf{ hyperfinite}, if for any $\eps>0$
there exists $K\geq 1$ such that
$\cP\subset \cH^d_{\eps,K}$.  
\vi
The significance of hyperfiniteness in algorithm theory
is highlighted by the following breakthrough result of Benjamini, Schramm and Shapira (Theorem 1.2 \cite{BSS}): every monotone hyperfinite property is testable
in constant-time (see \cite{GR} for property testing). Note that another proof of this statement
is given in \cite{HKNO}. It is important to note that
minor-closed families such as planar graphs, outer-planar graphs or bounded genus graphs are hyperfinite \cite{BSS}.
\vi
\section{Property A} \vi
In order to avoid confusion, in this section we use the phrase "graph class" instead of
"graph property", since we will talk about the notion of
Property A. \vi
First, let us formally define Property A.
Let $G\in\grd$ be a graph. Then, $\Prob(G)$ is the set
of all probability measures on the vertices of $G$.
If $f:V(G)\to \R$ and $g:V(G)\to \R$ are two real functions on the vertices then their $l_1$-distance is defined as $\| f-g\|_1:=\sum_{x\in V(G)} 
|f(x)-g(x)|$, also $\|f\|:=\sum_{x\in V(G)} |f(x)|\,.$
\begin{definition}[Property A]\label{prodef}
For $\eps>0$ and $r\geq 1$, a graph $G\in \grd$ is called \textbf{$(\eps,r)$-uniform} if there exists a probability measure valued function
$\tilde{f}:V(G)\to \Prob(G)$ such that
\begin{itemize}
\item for any adjacent pair of vertices $x,y\in V(G)$,
$\|\tilde{f}(x)-\tilde{f}(y)\|_1<\eps$,
\item for any $x\in V(G)$, we have that
$$\Supp(\tilde{f}(x))\subset B_r(x,G)\,,$$
\vi where $\Supp(\tilde{f}(x))$
denotes the set of vertices $z$
for which $\tilde{f}(x)(z)\neq 0.$
\end{itemize}
\noindent
We call a class of graphs $\cP$ $(\eps,r)$-uniform if all $G\in\cP$ are
$(\eps,r)$-uniform and we denote the class of all $(\eps,r)$-uniform graphs by $\cA^d_{\eps,r}$. 
A graph class $\cP\subset \grd$ is of \textbf{Property A} if for every $\eps>0$, there exists some $r\geq 1$ such that
$\cP\subset \cA^d_{\eps,r}$.
\end{definition}
\vi
Interestingly, Property A can be defined for a single
countable infinite graph or a finitely
generated group as well. Actually, the
notion of Property A has been introduced
for finitely generated groups by Guoliang Yu \cite{Yu} in the nineties with important
applications in algebraic topology and operator algebras \cite{NY}. It is not hard to see that
the set of paths $\{P_n\}^\infty_{n=1}$ forms a class of
Property~A, nevertheless later we will see that
the class of planar graphs is
of Property~A as well. 

\vi
As we will see in Section \ref{sia}, for finite graph classes Property A  is a strengthening of the notion
of hyperfiniteness. The proof of our main
theorem hinges on the fact that monotone hyperfinite graph classes are of
 Property A. 

\section{Property A and the Proof Labeling Schemes}
The sole goal of this section is to
prove the following proposition.
\begin{proposition}\label{relcheckame}
For any $0<\eps<\eps'<1$ and $r\geq 1$, $\cA^d_{\eps,r}$ can be verified by PLS in constant-time relative to $\cA^d_{\eps',r}$.
\end{proposition}
\proof A natural
approach for such a PLS is to label every vertex $x$ with its probability
distribution, described as a list of $\tilde{f}(x)(z)$, and let the vertices
check that the two conditions of Definition 4.1 hold.
There are two obstacles to this approach: the precise value of
$\tilde{f}(x)(z)$ might need a large number of bits to be encoded, and the
vertices do not agree on which vertex is "z" (remember that the vertices
do not have identifiers). The first lemma tackles the first obstacle
via discretization. For the second obstacle, we will use a coloring.
\begin{lemma} \label{csutlemma}
Let $G\in\grd$, $r \geq 1$, $x\in V(G)$, $f:B_r(x,G)\to \R$
be a nonnegative function such that $\sum_{y\in B_r(x,G)} f(y)=1$.
Let $\alpha> \frac{3}{\eps'-\eps}N^d_r$ be a positive integer (see Section \ref{masodiksection} for definition of $N^d_r$).
Then, there exists a function $g:B_r(x,G)\to \R$
such that
\begin{itemize}
\item $\sum_{y\in B_r(x,G)} g(y)=1$,
\item for any $y\in B_r(x,G)$, $g(y)=\frac{i}{\alpha}$, where
$0\leq i \leq \alpha$ is an integer,
\item $\sum_{y\in B_r(x,G)} |f(y)-g(y)|\leq \frac{\eps'-\eps}{3}.$
\end{itemize}
\end{lemma}
\proof
Let $g',g'':B_r(x,G)\to \R$ be defined in the following way.
$g'(y)=\frac{i}{\alpha}, g''(y)=\frac{i+1}{\alpha}$, where
$\frac{i}{\alpha}\leq f(y) \leq \frac{i+1}{\alpha}$. Then,
$\sum_{y\in B_r(x,G)} (f(y)-g'(y))< |B_r(x,G)|\frac{1}{\alpha}\leq\frac{\eps'-\eps}{3}$.
Similarly, $\sum_{y\in B_r(x,G)} (g"(y)-f(y)) < \frac{\eps'-\eps}{3}$.
\vi
Also, $$\sum_{y\in B_r(x,G)} g'(y)=\frac{k}{\alpha}\leq 1\,\,\,\mbox {and}\,\,
\sum_{y\in B_r(x,G)} g''(y)=\frac{l}{\alpha}\geq 1\,,$$
\vi
where $k$ and $l$ are integers such that $k\leq \alpha\leq l$. Note that for all $y\in B_r(x,G)$ we have $g''(y)-g'(y)=\frac{1}{\alpha}$, hence
$l-k=|B_r(x,G)|$. Pick a subset $S\subseteq B_r(x,G)$ such
that $|S|=\alpha-k$.
\vi
Let $g(y)=g'(y)$ if $y\notin S$ and $g(y)=g''(y)$ if $y\in S$.
Then,
$$\sum_{y\in B_r(x,G)}  g(y)=\frac{k}{\alpha}+\frac{\alpha-k}{\alpha}=1.$$
Also,
$$\sum_{y\in B_r(x,G)}  |f(y)-g(y)|\leq |B_r(x,G)| \frac{1}{\alpha}\leq \frac{\eps'-\eps}{3}\,.\,\,\qed$$
\vi
\begin{corollary} \label{corol1}
Let $0<\eps<\eps'<1$ and $G\in\cA^d_{\eps,r}$. Then, we have
a probability measure valued function
$\tilde{g}:V(G)\to \Prob(G)$ such that
\begin{itemize}
\item
for any $x\in V(G)$, 
$\tilde{g}(x)(z)=\frac{i}{\alpha}$, where
$0\leq i \leq \alpha$ is an integer
and $\alpha$ is the integer defined
by the previous lemma.
\item
for any adjacent pair
of vertices $x,y\in V(G)$,
$\|\tilde{g}(x)-\tilde{g}(y)\|_1<\epsilon'$,
\item
for any $x\in V(G)$, we have that
$\Supp(\tilde{g}(x))\subset B_r(x,G)\,.$
\end{itemize}
\end{corollary}
\proof 
Let $\tilde{f}:V(G)\to \Prob(G)$ the
probability measure valued function
in Definition \ref{prodef}.
For each $x\in V(G)$ we define
$g_x:V(G)\to \R$ in such a way that
$g_x$ satisfies the conditions of
the previous lemma with respect to the
function $\tilde{f}(x)$.
Now we define $\tilde{g}(x):=g_x$.\vi
Let $x\sim y\in V(G)$.
Then,
$$\|\tilde{g}(x)-\tilde{g}(y)\|_1\leq
\|\tilde{g}(x)-\tilde{f}(x)\|_1+
\|\tilde{f}(x)-\tilde{f}(y)\|_1+
\|\tilde{f}(y)-\tilde{g}(y)\|_1\leq$$
$$\leq 2\frac{\eps'-\eps}{3}+\eps<\eps'.\,\,\,\,\qed$$
\vi
Now we build the proof labeling scheme.
\begin{lemma} \label{color1}
Let $Q_1$ be a finite set 
such that $|Q_1|=N^d_{r}+1$.
Then, for any $G\in \grd$
there exists a coloring $S_1:V(G)\to Q_1$
such that if $d_G(x,y)\leq r$, then
$S_1(x)\neq S_1(y)$.
\end{lemma} 
\proof
Let $\hat{G}$ be the graph with vertex
set $V(G)$ such that $x,y\in V(G)$
are adjacent in $\hat{G}$ if and only
if $d_G(x,y)\leq r$. Then, the
degree of each vertex in $\hat{G}$
is at most $N^d_r$. Hence by
the classical Brooks' Theorem, there
is a proper coloring $S_1:V(G)\to Q_1$
of the graph $\hat{G}$. By definition,
$S_1$ defines a coloring of $G$
satisfying the condition of our lemma. \qed
\vi
Let $Q_2$ denote the finite
set $\{0,1,2,\dots,\alpha\}$, where
$\alpha$ is the integer in Lemma \ref{csutlemma}. Also,
let $Q_3$ be the set of
all maps $\phi:Q_1\to Q_2$.
Let $\pi_1:Q_1\times Q_3\to Q_1$ be
the first coordinate projection and
$\pi_2:Q_1\times Q_3\to Q_3$ be
the second coordinate projection.
For a $Q_1\times Q_3$-proof $T:V(G)\to Q_1\times Q_3$
let $T_1:V(G)\to Q_1$ be defined
as $\pi_1\circ T$ and let
$T_2:V(G)\to Q_3$ be defined
as $\pi_2\circ T$. 
\vi
We call $T:V(G)\to Q_1\times Q_3$
proper if $T_1(x)\neq T_1(y)$ provided
that $d_G(x,y)\leq r$.
If $T$ is proper then we can define
a function of two variables
$\tilde{T}:V(G)\times V(G)\to \R$
in the following way.
\begin{itemize}
\item If $d_G(x,y)\leq r$ then
$\tilde{T}(x,y)=\frac{i}{\alpha}$,
where $i=T_2(y)(T_1(x))$.
\item If $d_G(x,y)>r$ then
$\tilde{T}(x,y)=0$.
\end{itemize}
\vi
Properness is used to break possible symmetries. 
The verifier $\cV\subset B^d_{r+1,Q_1\times Q_3}$
is defined in the following way.
Let $N$ be a ball of radius at most $r+1$
and $C:V(N)\to Q_1\times Q_3$ be a $Q_1\times Q_3$-labeling on $N$.
Again, let $C_1=\pi_1\circ C$ and
$C_2=\pi_2\circ C$.
Then the $Q_1\times Q_3$-labeled
ball $(N,C)$ is in the verifier $\cV$ if
the following conditions are satisfied.
\begin{itemize}
\item {(Checking properness)}
If $y,z\in V(N)$ and $d_N(y,z)\leq r$
then $C_1(y)\neq C_2(z)$.
\item {(Checking probability)}
If $x$ is the center of $N$ then
$$\sum_{z\in B_r(x,N)} 
\frac{(C_2(z))(C_1(x))}{\alpha}=1\,.$$
\item {(Checking $l_1$-distance)}
If $x$ is the center of $N$ and
$x\sim y$, then
$$\sum_{z\in V(N)}
\left| \frac{ (C_2(z))(C_1(x))}{\alpha}
-\frac{ (C_2(z))(C_1(y))}{\alpha} \right|\leq \eps'.$$
\end{itemize}
\vi
Therefore, if $\cV$ accepts the
proof $T$, then
\begin{itemize}
\item $T$ is proper.
\item $\sum_{z\in V(G)} \tilde{f}(x)(z)=1$,
where $\tilde{f}(x)(z)=\tilde{T}(x,z)$.
\item For any adjacent pair $x\sim y\in V(G)$, $\|\tilde{f}(x)-\tilde{f}(y)\|\leq \eps'\,.$
\end{itemize}
\vi
Hence, if $\cV$ accepts $T$ then
$G\in \cA^d_{\eps',r}$.
\vi
In order to finish the
proof of the proposition, we
need to prove that if
$G\in \cA^d_{\eps,r}$
then there exists
a proof $T:V(G)\to Q_1\times Q_3$
that is accepted by the verifier $\cV$.
\vi
Let $\tilde{g}$ be the probability 
measure valued function defined
in Corollary \ref{corol1}.
Let $S_1:V(G)\to Q_1$ be the function
defined in Lemma \ref{color1}.
Finally, let $T:V(G)\to Q_1\times Q_3$
be defined in the following way.
\begin{itemize}
\item $T_1:=S_1$.
\item If $z\in V(G), q\in Q_1$ and there is
no $x\in B_r(z,G)$ such that $S_1(x)=q$,
then let $(T_2(z))(q)=0.$
\item If $z\in V(G), q\in Q_1$ and there exists
$x\in B_r(z,G)$ such that $S_1(x)=q$,
let $(T_2(z))(q)=i,$
where $\tilde{g}(x)(z)=\frac{i}{\alpha}.$
\end{itemize}
\vi
Then, $\cV$ accepts $T$.
Therefore, 
$$\cA^d_{\eps,r}\subset \cL_\cV\subset
\cA^d_{\eps',r}\,,$$
hence our proposition follows. \qed
\section{Strong hyperfiniteness}\vi
In this section we discuss some strengthenings of the notion of hyperfiniteness. \vi
\begin{definition} $G\in\grd$ is uniformly $(\eps,K)$-hyperfinite
if for all induced subgraph $F\subset G$, $F$ is $(\eps,K)$-hyperfinite
as well.
\end{definition}
\noindent
We say that a graph property $\cP\subset\grd$ is uniformly $(\eps,K)$-hyperfinite if all $G\in\cP$ are uniformly $(\eps,K)$-hyperfinite.
The set of all uniformly $(\eps,K)$-hyperfinite graphs will be
denoted by $\cU\cH^d_{\eps,K}$. We call a graph property $\cP$ \textbf{uniformly hyperfinite} if for any $\eps>0$, there exists $K\geq 1$ such that
$\cP\subset \cU\cH^d_{\eps,K}$.\vi
Monotone hyperfinite classes are, by definition, uniformly
hyperfinite, since they are closed to taking subgraphs. 
In our paper, we will use another strengthening of hyperfiniteness
introduced by Romero, Wrochna and \v{Z}ivn\'y \cite{RWZ} under the
name of \emph{fractional-cc-fragility}. First, we need a definition.
For a graph $G\in\grd$, we call $Y\subset V(G)$ a \emph{$K$-separator} if
by removing $Y$ (and all the adjacent edges) the components of the remaining graph
are of size at most $K$. We denote
the set of all $K$-separators of $G$ by
$\Sep(G,K)$.
\begin{definition} A graph $G\in\grd$ is strongly $(\eps,K)$-hyperfinite
if there exists a probability measure $\mu$ on $\Sep(G,K)$ such that for any $x\in V(G)$
$$\mu(Y\in \Sep(G,K)\,\mid\,x\in Y)<\eps\,.$$
\end{definition}
\noindent
We say that a graph class $\cP\subset\grd$ is strongly $(\eps,K)$-hyperfinite if all $G\in\cP$ are 
strongly $(\eps,K)$-hyperfinite.
The set of all strongly $(\eps,K)$-hyperfinite graphs will be
denoted by $\cS\cH^d_{\eps,K}$. We call a graph class $\cP$ \textbf{strongly hyperfinite} if for any $\eps>0$, there exists $K\geq 1$ such that
$\cP\subset \cS\cH^d_{\eps,K}$.\vi
It was first proved by Romero, Wrochna and \v{Z}ivn\'y (Theorem 1.5\,,\cite{RWZ}) that
monotone hyperfinite properties are strongly hyperfinite. The strong hyperfiniteness of monotone hyperfinite classes
plays an important role in the proof of Theorem \ref{main}.
Note that the author later proved in \cite{Elekuni} that the notions
of Property A, uniform hyperfiniteness and
strong hyperfiniteness in fact coincide.
\section{Property A implies hyperfiniteness}\label{amehyp}
\vi In this section we continue
the study of Property A and prove the
central technical proposition of our paper.
\vi
First, let us fix some notation, which will be used in the section.
Let $G\in\grd$ and $A\subset V(G)$. Then, $\partial_G(A)$ is the
set of vertices $x\in A$ such that there exists $y\notin A$, $x\sim y$.
Also,  $\partial^e_G(A)$ is the
set of edges $e=(x,y)$, where $x\in\partial_G(A)$ and $y\notin A$.
So, we have that
\begin{equation}\label{etrivi}
|\partial_G(A)|\leq |\partial^e_G(A)|.
\end{equation}
\begin{proposition}\label{zivn}
For any $\eps>0$ and $r\geq 1$,
$$\cA^d_{\eps,r}\subset \cH^d_{\frac{d^2\eps}{2}, N^d_{2r}},$$
\noindent
where $N^d_{2r}$ is defined in Section \ref{masodiksection}.
\end{proposition}
\proof
First, we need a technical lemma, which is very similar to Proposition 4.2 
in \cite{Tu}. Let $G\in\grd$ and $F\subset G$ be an induced
subgraph. We say that $F$ is $(\eps,r)$-uniform relative to $G$ 
if there exists a probability measure valued function $\tilde{f}: V(F)\to \Prob(F)$ such that
for any pair of adjacent vertices $x,y\in V(F)$,
\begin{equation}\label{kedd4}
\|\tilde{f}(x)-\tilde{f}(y)\|_1 \leq \eps
\end{equation}
\noindent
and for any $x\in V(F)$,
\begin{equation} \label{kedd5}
\Supp(\tilde{f}(x))\subset B_r(x,G).
\end{equation}
\begin{lemma}
If $G$ is $(\eps,r)$-uniform and $F\subset G$ is an induced subgraph,
then $F$ is $(\eps,2r)$-uniform relative to $G$.
\end{lemma}
\proof
For $x\in V(G)$, pick $\tau(x)\in V(F)$ in such a way
that $d_G(x,\tau(x))=d_G(x,F)$.
Let $g:V(G)\to\Prob(G)$ be a probability measure valued function
witnessing the fact that $G\in \cA^d_{\eps,r}$, that is,
\begin{itemize}  \label{e10171}
\item for any  adjacent pair of vertices $x,y\in V(G)$
\begin{equation}
\|\tilde{g}(x)-\tilde{g}(y)\|_1\leq \eps\,.
\end{equation} 
\item for any $x\in V(G)$
\begin{equation} \label{e10172}
\Supp(\tilde{g}(x))\subset B_r(x,G)\,.
\end{equation}
\end{itemize}
We define the function $\tilde{f}:V(F)\to \Prob(F)$ by $\tilde{f}(x)(z)=\sum_{t\in \tau^{-1}(z)} \tilde{g}(x)(t).$  Note that $\tau^{-1}(z)$ denotes the set of  vertices
mapped to $z$ by $\tau$.
Then by definition, $\Supp \tilde{f}(x)\subset V(F)$ and \\
for all $z\in V(F)$, $\tilde{f}(x)(z)\geq 0$. Also,
since $\cup_{z\in V(F)} \tau^{-1}(z)=V(G)$ we have that
$$\sum_{z\in V(F)} \tilde{f}(x)(z)= \sum_{t\in V(G)} \tilde{g}(x) (t)=1\,,$$
hence $\tilde{f}:V(F)\to \Prob(F)$. Also, if $x,y$ are adjacent vertices,
then
$$\|\tilde{f}(x)-\tilde{f}(y)\|_1 \leq \eps. $$
Indeed,  
$$\|\tilde{f}(x)-\tilde{f}(y)\|_1=\sum_{z\in V(F)} |\tilde{f}(x)(z)-\tilde{f}(y)(z)|\leq $$
$$\leq \sum_{z\in V(F)} |\sum_{t\in\tau^{-1}(z)} \tilde{g}(x)(t)-
\sum_{t\in\tau^{-1}(z)} \tilde{g}(y)(t)|\leq 
\sum_{z\in V(F)} \sum_{t\in\tau^{-1}(z)} |\tilde{g}(x)(t)-\tilde{g}(y)(t)|= $$
$$= \sum_{t\in V(G)} |\tilde{g}(x)(t)-\tilde{g}(y)(t)|=\|\tilde{g}(x)-\tilde{g}(y)\|_1\leq \eps.$$
\noindent
Also, 
\begin{equation} \label{e10173} \Supp (\tilde{f}(x))\subset B_{2r}(x,G).
\end{equation} Indeed,
if $\tilde{f}(x)(z)\neq 0 $, then there exists $t\in \tau^{-1}(z)$
such that $\tilde{g}(x)(t)\neq 0$. Hence by \eqref{e10172}, $d_G(t,x)\leq r$
and also, $d_G(t,z)\leq r$, since $d_G(t,z)\leq d_G(t,x)$ by the 
definition of $\tau$.  
That is, $d_G(x,z)\leq 2r$, so our lemma follows. \qed
\begin{lemma} Let $G\in \cA^d_{\eps,r}$ and let $F\subset G$ be
an induced subgraph. Then, there exists a non-empty subset
$L\subset V(F)$ such that $|\partial_F(L)|\leq \frac{d \eps}{2} |L|$
and $|L|\leq N^d_{2r} $.
\end{lemma}
\proof
Let $\tilde{f}:V(F)\to \Prob(F)$ be a probability measure valued function satisfying \eqref{kedd4} and \eqref{e10173}.
Such function exists by the previous lemma.
Then,
$$\sum_{x\in V(F)} \sum_{x\sim y} \|\tilde{f}(x)-\tilde{f}(y)\|_1\leq
\sum_{x\in V(F)} d \eps =$$
$$= \sum_
{x\in V(F)} d \eps \|\tilde{f}(x)\|_1\,.$$
\noindent
Hence,
$$\sum_{z\in V(F)} \sum_{x\in V(F)} \sum_{x\sim y}
|\tilde{f}(x)(z)-\tilde{f}(y)(z)|\leq d\eps \sum_{z\in V(F)}\sum_{x\in V(F)} 
\tilde{f}(x)(z)\,.$$
\noindent
Hence, there exists $z_0\in V(F)$ such that
$$\sum_{x\in V(F)} \sum_{x\sim y}
|\tilde{f}(x)(z_0)-\tilde{f}(y)(z_0)|\leq d \eps \sum_{x\in V(F)} 
\tilde{f}(x)(z_0)\,.$$
\noindent
We define the function $\zeta:V(F)\to [0,1]$ by
$\zeta(x)=\tilde{f}(x)(z_0)$, and we have that
\begin{equation}\label{kedd4uj}
\sum_{x\in V(F)} \sum_{x\sim y}
|\zeta(x)-\zeta(y)|\leq d \eps \sum_{x\in V(F)} 
\zeta(x)\,.
\end{equation}
\noindent
So far, we followed the proof of Proposition 3.2 in \cite{Brodzki}, however, in order to avoid some heavy machinery, we now choose a different
path.
Let us recall the area and coarea formulas (Lemma 3.6 and Lemma 3.7) from
\cite{KM}.
If $F\in\grd$ and $\zeta:V(F)\to [0,1]$, then we have the
following equations:
\begin{equation}
\label{kedd5uj}
\sum_{x\in V(F)} \sum_{x\sim y} |\zeta(x)-\zeta(y)|= 2 \int_0^1 
|\partial_F^e(\Omega_t(\zeta))|\,dt,
\end{equation}
\noindent
and
\begin{equation}
\label{kedd6uj}
\sum_{x\in V(F)}\zeta(x)= \int_0^1 |\Omega_t(\zeta)|\,dt,
\end{equation}
\noindent
where 
$$\Omega_t(\zeta)=\{x\in V(G)\,\mid\, \zeta(x)>t\}.$$
So, by \eqref{kedd4uj} 
$$ 2\int_0^1 |\partial_F^e(\Omega_t(\zeta))|\,dt\leq
d \eps \int_0^1 |\Omega_t(\zeta)|\,dt. $$
\noindent
Thus for some $t\geq 0$, we have
\begin{equation} \label{e10h}
|\partial_F^e(\Omega_t(\zeta))|\leq \frac{d \eps}{2} |\Omega_t(\zeta)|.
\end{equation}
\noindent
Now let $L=\Omega_t(\zeta)$. Then by \eqref{e10h} and \eqref{etrivi} we have that $|\partial_F(L)|\leq \frac{d \eps}{2} |L|$.
Also by definition, if $x\in L$, then $\tilde{f}(x)(z_0)>0$. Therefore,
$x\in B_{2r}(z_0,G)$. So, $0<|L|\leq N^d_{2r}$. Hence, 
our lemma follows. \qed \vi
Now we finish the proof of our proposition. Let $F_1=G$. Using the previous lemma, we choose
$L_1\subset V(F_1)$ to be a set of size at most $N^d_{2r}$ 
such that $|\partial_{F_1}(L_1)|\leq \frac{d \eps}{2} |L_1|$. Then, we
remove from $G$  all the edges outgoing from $L_1$. The number of such edges is
at most $d |\partial_{F_1}(L_1)|\leq \frac{d^2\eps}{2} |L_1|.$ Let $F_2$ be the subgraph of $G$ induced on $V(G)\backslash L_1$. Let $L_2\subset V(F_2)$ be a set of size at most $N^d_{2r}$ 
such that $|\partial_{F_2}(L_2)|\leq \frac{d \eps}{2}|L_2|$. Again, we
remove  from $G$ all the edges
outgoing from $L_2$. 
Inductively, we construct
disjoint components $L_1,L_2,\dots, L_n$ of size at most $N^d_{2r}$ such that $\cup^n_{i=1} L_i=V(G)$, by removing
at most $\frac{d^2\eps}{2} |V(G)|$ edges.
Hence, our proposition follows. \qed
\section{Strong hyperfiniteness implies Property A} \label{sia}
\begin{proposition}\label{hippo}
For any $\eps>0$ and $K\geq 1$, 
$$\cS\cH^d_{\eps,K}\subset \cA^d_{4\eps,K}\,.$$
\end{proposition}
\proof
Let $G\in \cS\cH^d_{\eps,K}$ and $\mu$ be a probability measure
on the set $\Sep(G,K)$ of $K$-separators of $G$ such that for any $x\in V(G)$ we
have
$$ \mu(Y\in \Sep (G,K)\,\mid\,x\in Y)\leq \eps\,.$$
\noindent
For a $K$-separator $Y$ and $x\in V(G)$, let the non-negative function
$f_{Y,x}:V(G)\to \R$ be defined in the following way.
\begin{itemize}
\item If $x\in Y$, let $f_{Y,x}(x)=\mu(Y)$ and
if $z\neq x$ let $f_{Y,x}(z)=0$.
\item If $x\notin Y$, then let $C_{Y,x}$ be the component
of $G\backslash Y$ containing the vertex $x$. If $z\in C_{Y,x}$ let $f_{Y,x}(z)=\frac{\mu(Y)}{|C_{Y,x}|}$. 
On the other hand, if  $z\notin C_{Y,x}$, let $f_{Y,x}(z)=0$.
\end{itemize}
\noindent
Then, for any $x\in V(G)$ and $Y\in \Sep(G,K)$ we have that
\begin{equation} \label{e11c}
\|f_{Y,x}\|_1=\sum_{z\in V(G)} f_{Y,x}(z)=\mu(Y)\,.
\end{equation}
\vi
Hence, 
$$\tilde{f}(x)=\sum_{Y\in \Sep(G,K)} f_{Y,x}\,$$
\vi defines a probability measure valued function $\tilde{f}:V(G)\to \Prob(G)$.
\noindent
Then, by the definition of $K$-separators,
 for any $x\in V(G)$, 
 \begin{equation} \label{e11a} \Supp(\tilde{f}(x))\subset B_K(x,G)\,.\end{equation}
 Also, for any pair of adjacent vertices $x\sim y$, we have that
 \begin{equation} \label{e11b}
\|\tilde{f}(x)-\tilde{f}(y)\|_1 \leq 4\eps\,.
\end{equation}
\noindent
Indeed, if $x,y\notin Y$ and $x\sim y$, then by definition, $f_{Y,x}=f_{Y,y}$. 
If $x\in Y$ or $y\in Y$, then by \eqref{e11c},
$$\|f_{Y,x}-f_{Y,y}\|_1\leq \|f_{Y,x}\|_1+ \|f_{Y,y}\|_1\leq 2 \mu(Y)\,.$$
Therefore,
$$\|\tilde{f}(x)-\tilde{f}(y)\|_1 \leq \sum_{Y\in \Sep(G,K)} \|f_{Y,x}-f_{Y,y}\|_1 \leq $$
$$ \leq \sum_{Y\in \Sep(G,K), x\in Y} 2\mu(Y) +  \sum_{Y\in \Sep(G,K), y\in Y} 2\mu(Y) + \sum_{Y\in \Sep(G,K), x,y\notin Y} \|f_{Y,x}-f_{Y,y}\|_1 \leq
 4\eps\,,$$
 \vi
since if $x\sim y$ and $x,y\notin Y$, then $f_{Y,x}=f_{Y,y}$.
Hence by \eqref{e11a} and \eqref{e11b}, $G\in \cA^d_{4\eps,K}$. \qed
\begin{lemma}\label{corolfontos}
For any $\eps>0$, there exists $r>0$ such that
$\cP\subset \cA^d_{\eps,r}$, where $\cP\subset \grd$ is a monotone hyperfinite class of graphs.
\end{lemma}
\proof
By Theorem 1.6 of \cite{RWZ}, for any $\eps>0$, there exists $K>0$
such that $\cP\subset \cS\cH^d_{\eps,K}$. Hence our lemma follows, from
Proposition \ref{hippo}.\qed

\section{The Proof of the Main Theorem}
In this section, we prove Theorem \ref{main}.
Let  $\cP\subset \grd$ be a monotone hyperfinite property
that is closed under taking disjoint unions and let $\eps>0$.
We need to prove that there exists a verifier $\cV_\eps$ such that (using
the notation of Section \ref{masodiksection})
\begin{equation} \label{kulcs}
\cP\subset \cL_{\cV_\eps} \subset \cP_\eps\,.
\end{equation}
\noindent
For $r>0$, let $L\cP_r\subset \grd$ denote the class of \emph{$r$-locally $\cP$ graphs.}  That is, $G\in L\cP_r$ if
for all $x\in V(G)$, $B_r(x,G)\in \cP$.
So, by definition, there exists a verifier
$\cB_K\subset B^d_{K,Q}$ such that $L\cP_K=\cL_{\cB_K}$ and
the finite set $Q$ is empty.
\begin{lemma}\label{mostl1}
For any $\eps>0$, there exists $M_\eps>0$ such that
if $K\geq M_\eps$, then
$$\cP\subset \cH^d_{\eps,K}\cap L\cP_{K}\subset
\cP_{\eps}\,.$$
\end{lemma}
\proof
Pick $M_\eps>0$ such that $\cP\subset \cH^d_{\eps,M_\eps}$, let 
$K\geq M_\eps$. Then, we have $\cP\subset \cH^d_{\eps,K}\cap L\cP_{K}$. 
Now, let $G\in \cH^d_{\eps,K}\cap L\cP_{K}$. So, we can remove at most
$\eps |V(G)|$ edges from $G$ to obtain a graph having components
of size at most $K$ and by monotonicity, all those components are in $\cP$. Also, by our assumption about $\cP$,
the disjoint unions of the components are in $\cP$ as well.
Hence, $G\in \cP_\eps$. That is, $\cH^d_{\eps,K}\cap L\cP_{K}\subset
\cP_{\eps}\,.$ \qed
\begin{lemma}\label{mostl2}
For any $\eps>0$, there exists $N_\eps>0$ and a verifier
$\cC_\eps$ such that if $K\geq N_\eps$, then
$$\cP\subset \cL_{\cC_\eps} \subset \cH^d_{\eps,K}\,.$$
\end{lemma}
\proof 
Let $r\geq 1$ be an integer such that $\cP\subset \cA^d_{\frac{\eps}{d^2},r}$.
Such $r$ exists by Lemma \ref{corolfontos}. 
Let $N_\eps=N^d_{2r}$. By Proposition \ref{relcheckame},
there exists a verifier $\cC_\eps$ such that
$$\cA^d_{\frac{\eps}{d^2},r}\subset \cL_{\cC_\eps}
\subset \cA^d_{\frac{2\eps}{d^2},r}\,.$$
\vi
By Proposition \ref{zivn}, 
$$\cA^d_{\frac{2\eps}{d^2},r}\subset \cH^d_{\eps,N^d_{2r}}=
\cH^d_{\eps, N_\eps}\,.$$
\vi
Therefore, 
$$\cP\subset \cL_{\cC_\eps} \subset \cH^d_{\eps,K}\,.\,\,\,\qed$$
\vi
\noindent
Now we finish the proof of Theorem \ref{main}.
Let $K_\eps=\max (N_\eps, M_\eps)$ and $\cD_\eps=\cB_{K_\eps}$,
then by our previous lemmas,
$$\cP\subset \cL_{\cC_\eps}\cap \cL_{\cD_\eps}\subset \cP_\eps\,.$$
\noindent
Hence by Lemma \ref{szorzat}, we have a verifier $\cV_\eps$
satisfying \eqref{kulcs}. \qed

\end{document}